\documentclass[12pt,twoside]{amsart}
\usepackage{amsmath}
\usepackage{amsthm}

\usepackage{amsfonts}
\usepackage{amssymb}
\usepackage{latexsym}
\usepackage[all]{xy}

\date{}
\allowdisplaybreaks[4] \footskip=15pt
\renewcommand{\uppercasenonmath}[1]{}

\numberwithin{equation}{section} \theoremstyle{plain}
\newtheorem*{thm*}{Main Theorem}
\newtheorem{thm}{Theorem}[section]
\newtheorem{cor}[thm]{Corollary}
\newtheorem*{cor*}{Corollary}
\newtheorem{lem}[thm]{Lemma}
\newtheorem*{lem*}{Lemma}
\newtheorem{prop}[thm]{Proposition}
\newtheorem*{prop*}{Proposition}

\newtheorem*{rem*}{Remark}
\newtheorem{exa}[thm]{Example}
\newtheorem*{exa*}{Example}

\newtheorem*{df*}{Definition}

\newtheorem*{conj*}{Conjecture}
\newtheorem*{ack*}{ACKNOWLEDGEMENTS}

\newcommand{\pf}{\noindent\begin {proof}}
\newcommand{\epf}{\end{proof}}

\begin{document}
\title{A Note on Quasi-Frobenius Rings}
\author{ Liang Shen and Jianlong Chen}

\address{Department of Mathematics, Southeast University, Nanjing, 210096,
P.R.China}

\address{Department of Mathematics, Southeast University, Nanjing, 210096,
P.R.China}

%\date{December, 2004}
\maketitle \baselineskip=20pt
\begin{abstract}
The Faith-Menal conjecture  says that every strongly right $Johns$
ring is $QF$.
 The conjecture is also equivalent to say every right noetherian left $FP$-injective ring
 is $QF$. In this short article, we show that the conjecture is true under the
 condition( a proper generalization of left $CS$ condition)
  that every nonzero  complement left ideal is not small( a left ideal $I$ is called small
  if for every left ideal $K$, $K$+$I$=$R$ implies $K$=$R$). It is also proved that (1) $R$
  is $QF$ if and only if $R$ is a left and right mininjective ring with $ACC$ on right
annihilators in
  which $S_{r}\subseteq ^{ess}R_{R}$; (2) $R$ is $QF$ if and only
  if $R$ is a right simple injective ring with $ACC$ on right annihilators in
  which $S_{r}\subseteq ^{ess}R_{R}$. Several known results on $QF$ rings are obtained as
corollaries.
\end{abstract}

\bigskip
\section { \bf INTRODUCTION}
\bigskip

Throughout this paper rings are associative with identity, a
regular ring means a von Neumann regular ring. We write $J=J(R),
Z_l, Z_r, S_l$ and $S_r$ for the Jacobson radical, the left
singular ideal, the right singular ideal, the left socle and the
right socle of $R$ respectively. We use $X\subset Y$ to mean that
the inclusion is proper. For a subset $X$ of a ring $R$, the left
annihilator of $X$ in $R$ is ${\bf l}(X)=\{r\in R: rx=0$ for all
$x\in X\}$. For any $a\in R$, we write ${\bf l}(a)$ for ${\bf
l}(\{a\})$. Right annihilators are defined analogously.

  A ring $R$ is called  left $n$-injective ( left mininjective) if, for any $n$-generated (minimal) left ideal $I$ of
$R$, every $R$-homomorphism from $I$ to $R$ extends to an
$R$-homomorphism from $R$ to $R$. Left {\it1}-injective rings are
called left $P$-injective, or equivalently, ${\bf rl}(a)$=$aR$ for
each $a\in R$. A ring $R$ is said to be right simple injective if
every homomorphism from a right ideal of $R$ to $R$ with simple
image can be given by left multiplication by an element of $R$. It
is clear that right simple injective rings imply right
mininjective rings. A ring $R$ is called left {\it FP}-injective
if, for any free left $R$-module $F$ and any finitely generated
$R$-submodule $N$ of $F$, every $R$-homomorphism $f: N\rightarrow
R$ can be extended to an $R$-homomorphism $g: F\rightarrow R$. Or
equivalently, the matrix ring M$_{n\times n}$($R$) is left
$P$-injective for each $n \geq 1$(see \cite[Theorem 5.41]{NY03}).
\\\indent It is well known that a ring $R$ is quasi-Frobenius (or $QF$) if and
only if $R$ is left or right noetherian and left or right
self-injective (see \cite[Theorem 30.10]{AF92}). There are three
open conjectures on $QF$ rings, which have attracted many people
(such as Faith,  Nicholson,  Yousif and so on) to work on them.
The three conjectures were deeply discussed in \cite{ FH02,NY03}.
One of the three conjectures is the Faith-Menal conjecture. It was
raised by  Faith and  Menal in \cite{FM94}. It says that every
strongly right Johns ring is $QF$. Recall that a ring $R$ is
called right Johns if $R$ is right noetherian and every right
ideal is an annihilator. Right Johns rings were characterized by
Johns in \cite{J77}, but he used a false result of Kurshan
\cite[Theorem 3.3]{K70} to show that right Johns rings are right
artinian. In \cite{FM92}, Faith and Menal gave a counter example
to show that right Johns rings need not be right artinian. Later
(see \cite{FM94}) they defined strongly right Johns ring(the
matrix ring M$_{n\times n}$($R$) is right Johns for all $n \geq
1$) and characterized such rings as right noetherian and left
$FP$-injective rings . But they didn't know whether a strongly
right Johns ring is $QF$. One key to prove the conjecture is to
show such ring is semilocal, or only need to prove that $R/J$ is
regular.
\\\indent Let $S$ be a submodule of an
$R$-module $_{R}M$. A submodule $C\subseteq M$ is said to be a
complement to $S$ (in $M$) if $C$ is maximal with respect to the
property that $C\cap S$=0. We say that a submodule $C\subseteq M$
is a complement in $M$ if there exists a submodule $S\subseteq M$
such that $C$ is a complement to $S$ in $M$. We prove that every
right noetherian and left $2$-injective ring is $QF$ if and only
if every non-zero complement left ideal is not small(see Theorem
2.4). Since left $FP$-injectivity implies $2$-injectivity by
definition, the Faith-Menal conjecture is true under the condition
that every non-zero complement left ideal is not small. Using a
similar way to prove the regularity of $R/J$, we also prove that
if $R$ is a right mininjective ring with $ACC$ on right
annihilators in which $S_{r}\subseteq ^{ess}R_{R}$, then $R$ is
semiprimary (see Lemma 2.11). By this useful lemma, we have (1) if
$R$ is a left and right mininjective ring with $ACC$ on right
annihilators in which $S_{r}\subseteq ^{ess}R_{R}$, then $R$ is
quasi-Frobenius(see Theorem 2.13); (2) if $R$ is a right simple
injective ring with $ACC$ on right annihilators in which
$S_{r}\subseteq ^{ess}R_{R}$, then $R$ is quasi-Frobenius(see
Theorem 2.15). Therefore we remove some conditions of four
theorems of Nicholson and Yousif (\cite[Theorem 3.31, Theorem
6.44, Theorem 8.4, Theorem 8.5 ]{NY03}). Several known results on
$QF$ rings are obtained as corollaries.

\bigskip
\section { \bf  RESULTS}
\bigskip

\begin{prop} \label{prop: 2.1} Let $R$ be a right noetherian and
left P-injective ring, then $J=$$Z_{l}$ is a right annihilator,
nilpotent and {\bf l}$(J)$ is essential both as a left and a right
ideal of $R$.
\end{prop}
\begin{proof}
By (\cite [Theorem 5.14]{NY03}), $J$=$Z_{l}$.  $J$ is nilpotent
and {\bf l}$(J)$ is essential both as a left and a right ideal of
$R$ followed by \cite[Theorem 2.7]{GG98}. Since {\bf l}$(J)$ is an
 essential left ideal, {\bf r}{\bf l}$(J)\subseteq$ $Z_{l}=J$, which
implies $J={\bf r}{\bf l}(J)$.
\end{proof}
\medskip
An element $a \in R$ is called regular if there exists an element
$b$ in $R$ such that $a=aba$.
\begin{lem}\label{lem: 2.2}
Let b $=$ $a-aca$ such that a, c $\in R$. If b is a regular
element, so is a.
\end{lem}
\begin{proof}
Since $b$ is a regular element, there exists $d \in R$ such that
$a-aca=b=bdb=(a-aca)d(a-aca)$ which shows that
$a=a[c+(1-ca)d(1-ac)]a$ is a regular element.
\end{proof}

\begin{lem}\label{lem: 2.3}
Let $R$ be a right noetherian and left P-injective ring such that
each non-zero complement left ideal is not small, then $R$ is
right artinian.
\end{lem}
\begin{proof}
First we prove that $\overline{R}=R/J$ is a regular ring. Assume a
$\notin J$,  since $J={\bf r}{\bf l}(J)$$=$$Z_{l}$ by Proposition
2.1, there exists a non-zero left ideal $I$, which is a complement
to {\bf l}($a$). We claim that there must exist $b\in$ $I$ such
that $ba \notin J$. If not, then $Ia\subseteq J$, which implies
that ${\bf l}(J)Ia=0$. Since ${\bf l}(a) \cap I=0$, {\bf
l}($J$)$I\subseteq {\bf l}(a)\cap I=0$. Thus $I \subseteq {\bf
r}{\bf l}(J)$=$J$.  So $I$ is a small left ideal, a contradiction.
As $R$ is left $P$-injective,
 every homomorphism from $Rba$ to $R$ is a right
multiplication. Define $f(rba)=rb$, then $f$ is well-defined and
there exists $0\ne c\in R$ such that $f(rba)=rbac$.
 Therefore $b=bac$, which implies $\overline b\in {\rm\bf l}(\overline a-\overline {aca})$. Since
 $\overline{ba}\neq\overline{0}$, {\bf l}($\overline{a}$)
  is properly contained in {\bf l}($\overline{a-aca})$. If $a-aca\in
 J$, then $\overline{a}$ is a regular element in
 $\overline{R}$. If not, let $a_{1}=a-aca$, then {\bf l}($a_{1}$) is not
 essential in $_{R}R$. In the same way we get
 $a_{2}=a_{1}-a_{1}c_{1}a_{1}$ for some $c_{1}\in R$ and
 {\bf l}($\overline{a_{1}}$) is properly contained in {\bf
 l}($\overline{a_{2}}$). If $a_{2}\in J$, then $\overline a_{1}$
 is a regular element in $\overline R$. Thus by Lemma 2.2,
 $\overline a$ is a regular element in $J$. If $a_{2}\notin J$,
 then we have $a_{3}=a_{2}-a_{2}c_{2}a_{2}$ for some $c_{2}\in R$.
 Therefore we have  such $a_{k}\in R$ step by step, $k=1,2,\cdots$. And we
 are going to be stopped after finite steps. If not,
 we get a chain of left annihilators in $\overline R$: {\bf l}($\overline a_{1}$)$\subset${\bf l}($\overline a_{2}$)
 $\subset\cdots$. Since $R$ is right noetherian and $J(\overline R)$= $\overline 0$,
  $\overline R$ is a semiprime and right Goldie ring. Thus $\overline R$ satisfies $ACC$ on left
 annihilators by \cite[Lemma 5.8]{GW89}, a contradiction. Then there must exist some
 positive integer $m$ such that $a_{m}\in J$ and $a_{k}=a_{k-1}-a_{k-1}c_{k-1}a_{k-1}$ for some $c_{k-1}\in R, k=2,3,\cdots,m$.
  Hence it is clear  that $\overline a$ is a regular element in $\overline
 R$ by Lemma 2.2.
 Since $\overline a$ is arbitrary, $\overline{R}$ is a regular ring. Then $\overline R$ is
 semisimple because $\overline R$ is right noetherian. Thus $R$ is semiprimary for $J$ is
 nilpotent by Proposition 2.1. So $R$ is right artinian.
\end{proof}
 A ring $R$ is called left $CS$ (left min-$CS$) if every
left ideal (minimal left ideal) is essential in a direct summand
of $_{R}R$. $R$ is said to be left $C2$ if every left ideal
isomorphic to a direct summand of $_{R}R$ is itself a direct
summand of $_{R}R$. A left $CS$ left $C2$ ring is called left
continuous. Every left self-injective ring is left $P$-injective
left $CS$; and every left $P$-injective left $CS$-ring is left
continuous because every left $P$-injective ring is left $C2$ by
\cite[Proposition 5.10]{NY03}. Left $CS$ condition is also
equivalent to say that every complement left ideal is a direct
summand. Then it is clear that left $CS$ condition implies that
every non-zero complement left ideal is not small, but the
converse is not true (see Example 2.8).

\begin{thm}\label{thm: 2.4}
R is QF if and only if R is right noetherian, left 2-injective and
every nonzero  complement left  ideal is not small.
\end{thm}
\begin{proof}
{\rm(i)} $``\Longrightarrow"$ Since $R$ is left self-injective, it
is left $CS$.  Therefore every nonzero complement left ideal $I$
is generated by a nonzero idempotent of $R$, which implies $I$ is
not small. The left is obvious.
\\\indent\indent$~$ {\rm(ii)}$``\Longleftarrow"$ By Lemma 2.3,  $R$ is right
artinian, then $R$ has $ACC$ on left annihilators. Thus $R$ is
$QF$ by \cite[Corollary 3]{R75}.
\end{proof}
\begin{cor}\label{cor: 2.5}
If $R$ is right noetherian and left FP-injective, then R is QF if
and only if  every nonzero complement left ideal is not small.
\end{cor}

\begin{thm}\label{thm: 2.6}
If $R$ is right noetherian, left P-injective, every nonzero
complement left ideal is not small and left min-CS, then $R$ is
QF.
\end{thm}
\begin{proof}
By Lemma 2.3, $R$ is right artinian, then $R$ is a left $GPF$
ring(left $P$-injective, semiperfect, and $S_{l}\subseteq
^{ess}$$_{R}R$). So $R$ is left Kasch and $S_{r}=S_{l}=S$ by
\cite[Theorem 5.31]{NY03}, then  $Soc(Re)$ is simple for each
local idempotent $e\in R$(see \cite[Lemma 4.5]{NY03}). Thus
$(eR/eJ)^{\ast}$$\cong l(J)e=S_re=Soc(Re)$ is simple for every
local idempotent $e\in R$. Since $R$ is semiperfect, each simple
right $R$-module is isomorphic to $eR/eJ$ for some local
idempotent $e\in R$ by \cite [Theorem 27.10]{AF92}. Hence $R$ is
right mininjective by \cite[Theorem 2.29]{NY03}. From above, $R$
is a two-sided mininjective and right artinian ring, then it is
$QF$ by \cite[Theorem 3.31]{NY03}.
\end{proof}
\begin{cor}\label{cor: 2.7}
\cite[Theorem 2.21]{CL04} If $R$ is right noetherian left CS and
left P-injective, then $R$ is QF.
\end{cor}
\begin{exa}\label{exa: 2.8}
There is a left min-CS ring $S$ satisfying every nonzero
complement left ideal is not small. But it is not left
$CS$.\end{exa}
\begin{proof} Let $k$ be a division
ring and $V_{k}$ be a right $k$-vector space of infinite
dimension. Take $R$=End($V_{k}$), then it is well-known that $R$
is regular but not left self-injective. Let $S$=M$_{2\times2}(R)$,
then $S$ is also regular, which implies $S$ is left $P$-injective
and $J(S)$ is zero. Thus $S$ is left C2, every non-zero complement
left ideal of $S$ is not small and every minimal left ideal is a
direct summand of $S$. But $S$ can not be left $CS$. For if $S$ is
left $CS$, then $S$ is left continuous. Hence $R$ is left
self-injective by \cite[Theorem 1.35]{NY03}, a contradiction.
\end{proof}
\begin{cor}\label{cor: 2.9}
\cite[Theorem 3.2]{NY98} The following condition on a ring R are
equivalent. \\(1) R is quasi-Frobenius. \\(2) R is a right Johns,
left CS-ring.
\end{cor}
\begin{lem}\label{lem: 2.10}
\cite[Lemma 2]{L69} Let R satisfy ACC on right annihilators, and
suppose {\bf l}(S) is a two-sided ideal of R. Then R/{\bf l}(S)
has ACC on right annihilators.
\end{lem}

 \indent By using a similar proof in Lemma 2.3 to show that
$R/J$ is regular, we get a useful lemma below.
\begin{lem}\label{lem: 2.11}
If $R$ is a right mininjective ring with $ACC$ on right
annihilators in which $S_{r}\subseteq ^{ess}R_{R}$, then $R$ is
semiprimary.
\end{lem}
\begin{proof}
First we prove $J={\bf l}(S_{r})=Z_{r}$. Since $R$ is right
mininjective, $S_{r}\subseteq S_{l}$ by \cite [Theorem 2.21]
{NY03}. So $J\subseteq {\bf l}(S_{r})$. ${\bf l}(S_{r})\subseteq
Z_{r}$ is followed by that $S_{r}\subseteq ^{ess}R_{R}$. Note that
$R$ has $ACC$ on right annihilators, hence $Z_{r}$ is nilpotent by
\cite[Theorem 7.15]{L98}. Therefore $Z_{r}\subseteq J$. Thus
$J={\bf l}(S_{r})=Z_{r}$ is nilpotent. Second we claim that
$\overline{R}$=$R/J$ is regular. If $\overline a\neq\overline 0$,
since $J=Z_{r}$, there exists a nonzero right ideal $I$ such that
{\bf r}($a$)$\cap I$=0. Note that $S_{r}\subseteq ^{ess}R_{R}$.
There exists a minimal right ideal $bR\subseteq I$. So {\bf
r}($a$)$\cap bR$=0. Define $f$: $abR\rightarrow bR$ by
$f(abr)=br$, $\forall r\in R$. It is clear that $f$ is
well-defined. Since $R$ is right mininjective and $abR$ is also a
minimal right ideal, there exists $c\in R$ such that $f(abr)=cabr$
for every $r\in R$. Thus $b=cab$. Therefore $b\in$ ${\bf
r}(a-aca)$$\backslash {\bf r}(a)$. If $a-aca\in J$, $\overline a$
is a regular element. If not, let $a_{1}=a-aca$. In the same way
we get $a_{2}=a_{1}-a_{1}c_{1}a_{1}$ for some $c_{1}\in R$ and
${\bf r}(a_{1})$ is properly contained in ${\bf r}(a_{2})$. Since
$R$ has $ACC$ on right annihilators, $\overline R$ is a regular
ring by a similar proof in Lemma 2.3. Next we show that $R$ is
semilocal. Since $J={\bf l}(S_{r})$ is a left annihilator, the
ring $\overline R$ inherits the $ACC$ on right annihilators by
Lemma 2.10. Thus $R/J$ is semisimple. Therefore $R$ is semilocal.
From above $R$ is semilocal and $J$ is nilpotent. So $R$ is
semiprimary.
\end{proof}
\indent Now we have a look at a theorem of Nicholson and Yousif,
which is a sharp improvement on Ikeda's theorem.
\begin{lem}\label{lem: 2.12}
\cite[Theorem 3.31] {NY03}
 Suppose that R is a semilocal, left and
right mininjective ring with $ACC$ on right annihilators in which
$S_{r}\subseteq ^{ess}R_{R}$. Then R is quasi-Frobenius.
\end{lem}
\indent By Lemma 2.11 and Lemma 2.12, the semilocal condition can
be removed and the ring satisfies the assumptions of \cite[Theorem
8.5]{NY03} is just $QF$. Thus we have
\begin{thm}\label{thm: 2.13}
If R is a left and right mininjective ring with $ACC$ on right
annihilators in which $S_{r}\subseteq ^{ess}R_{R}$, then R is
quasi-Frobenius.
\end{thm}

\noindent {\bf Remark.} None of the conditions above can be removed.\\
(1) $S_{r}\subseteq ^{ess}R_{R}$ is needed. For example, the ring
of integers $\mathbb{Z}$ is a two-sided mininjective and
noetherian ring with $S_{r}=0$. But it is not a $QF$ ring.\\
(2) Right(or left) mininjective can't be removed. For Bj\"{o}rk
Example(see \cite[Example 2.5]{NY03}) is a right mininjective and
left artinian ring. So it has $ACC$ on both right and left
annihilators.  And the two socles are both essential. But it is
not left mininjective.\\
(3) The chain condition is essential. As there exists a ring which
is left $PF$ (semiperfect, left self-injective and $S_{l}\subseteq
^{ess}$$_{R}R$) but not right $PF$(see \cite {DM86}). Such ring
satisfies the conditions except the chain condition in Theorem
2.13. But it is not $QF$.

\begin{cor}\label{cor: 2.14}
\cite[Theorem 8.4]{NY03} The following are equivalent for a ring
R:
\\(1) R is right and left mininjective, right noetherian, and $S_{r}\subseteq
^{ess}R_{R}$.
\\(2) R is right and left mininjective, right finitely
cogenerated, with ACC on right annihilators.
\\(3) R is quasi-Frobenius.
\end{cor}

\begin{thm}\label{thm: 2.15}
If R is  a right simple injective ring with $ACC$ on right
annihilators in which $S_{r}\subseteq ^{ess}R_{R}$. Then R is
quasi-Frobenius.
\end{thm}
\begin{proof}
By Lemma 2.11, $R$ is semiprimary. Then $R$ is left $P$-injective
by \cite[Theorem 6.16]{NY03}, which implies that $R$ is left
mininjective. Thus $R$ is $QF$ by Theorem 2.13.
\end{proof}
\indent Recall that a ring $R$ is called a right Goldie ring if it
has the $ACC$ on right annihilators and $R_{R}$ is finite
dimensional.
\begin{cor}\label{cor: 2.16}
\cite[Theorem 6.44] {NY03} The following conditions are equivalent
for a ring R: \\(1) R is quasi-Frobenius. \\(2) R is a right
simple injective, right noetherian ring with $S_{r}\subseteq
^{ess}R_{R}$. \\(3) R is a right simple injective, right Goldie
ring with $S_{r}\subseteq ^{ess}R_{R}$.
\end{cor}

\end{document}